\documentclass[a4paper]{elsarticle}
\usepackage{lipsum}
\makeatletter
\def\ps@pprintTitle{%
 \let\@oddhead\@empty
 \let\@evenhead\@empty
 \def\@oddfoot{}%
 \let\@evenfoot\@oddfoot}
\makeatother
\usepackage[left=3cm,right=3cm,top=2cm,bottom=2cm,includeheadfoot]{geometry}
\usepackage{amsmath}
\usepackage{amsfonts}
\usepackage{amssymb}
\usepackage{amsfonts}
\usepackage{ MnSymbol}
\usepackage{amsmath}
\usepackage{cancel}
\usepackage{graphicx}
\usepackage[numbers]{natbib}

\usepackage[ansinew]{inputenc}
\usepackage{amsthm}
\usepackage[arrow, matrix, curve]{xy}
\newtheoremstyle{style}   
 {0.5cm}                 
 {0.5cm}                 
  {}                         
  {}                         
  {\normalfont\bfseries}  
  {\normalfont }{   }
  {}
\newtheorem{thm}{Theorem}
\newtheorem{satz}[thm]{Theorem}
\newtheorem{Proposition}[thm]{Proposition}
\newtheorem{Lemma}[thm]{Lemma}

\theoremstyle{remark}
\newtheorem{bemerkung}{Remark}
\newtheorem{beispiel}{Example}

\newcommand{\ph}{\varphi}
\newcommand{\R}{\mathbb{R}}

\newcommand{\De}{\Delta_{p,q}}

\allowdisplaybreaks[1]
\begin{document}
\begin{frontmatter}

\title{\textbf{{Reducible conformal holonomy in any metric signature and application to twistor spinors in low dimension}}}


\author{Andree Lischewski}

\address{Humboldt-Universit{\"a}t zu Berlin, Institut f{\"u}r Mathematik, Rudower Chausse 25, 12489 Berlin, Germany}
\date{today}
\begin{abstract}
{\small{We prove that given a pseudo-Riemannian conformal structure whose conformal holonomy representation fixes a totally lightlike subspace of arbitrary dimension, there is, wrt. a local metric in the conformal class defined off a singular set, a parallel, totally lightlike distribution on the tangent bundle which contains the image of the Ricci-tensor. This generalizes results obtained for invariant lightlike lines and planes and closes a gap in the understanding of the geometric meaning of reducibly acting conformal holonomy groups. We show how this result naturally applies to the classification of geometries admitting twistor spinors in some low-dimensional split signatures when they are described using conformal spin tractor calculus. Together with already known results about generic distributions in dimensions 5 and 6 we obtain a complete geometric description of local geometries admitting real twistor spinors in signatures $(3,2)$ and $(3,3)$. In contrast to the generic case where generic geometric distributions play an important role, the underlying geometries in the non-generic case without zeroes turn out to admit integrable distributions.}}
\end{abstract}

\ead{lischews@math.hu-berlin.de}

\end{frontmatter}

\section{Introduction}

Modelling a pseudo-Riemannian conformal structure $(M,c=[g])$ of signature $(p,q)$,  where $n=p+q \geq 3$, as a parabolic Cartan geometry $(\mathcal{P}^1, \omega^{nc})$ of type $(G=O(p+1,q+1),P)$, where $P \subset G$ is the stabilizer of some isotropic ray, in the sense of \cite{sharp,cs,baju} leads to a well-defined algebraic conformal invariant, being the conformal holonomy group $Hol(M,c)$. As no canonical connection for $(M,c)$ can be defined on a reduction of the frame bundle of $M$, the Cartan geometry in question arises via a procedure called the first prolongation of a conformal structure, which naturally identifies $Hol(M,c) \cong Hol(\omega^{nc})$ with a (class of conjugated) subgroup of $O(p+1,q+1)$. Conformal holonomy groups turn out to be interesting objects in their own right and as in the metric case one is particularly interested in the relation between algebraic properties of $Hol(M,c)$ and underlying geometric structures in the conformal class, cf. \cite{arm,leihabil,baju,ln,lst,alt}, for instance.
As the construction via Cartan geometries differs from the definition of $Hol(M,g) \subset O(p,q)$ for the metric case, one is met with new challenges and features when carrying out the above programme:\\
A complete classification of irreducibly acting conformal holonomy groups is hindered by the fact that there is no useful algebraic criterion for $\mathfrak{hol}(M,c)$, such as being a Berger algebra for metric holonomy groups, which would reduce the problem to a finite classification list. However, \cite{alt} classifies all irreducibly acting conformal holonomy groups which at the same time also act transitively on the M{\"o}bius sphere in signature $(p,q)$, being the projectivized null-cone in $\R^{p+1,q+1}$. Among other groups, one finds in this list special unitary conformal holonomy, which means that there is locally a Fefferman spin space over a strictly pseudoconvex spin manifold in the conformal class (cf. \cite{baju,leihabil}). Conformal structures with $Hol(M,c) \subset Sp(k+1,m+1)$ were studied in detail in \cite{altdr}. The models of such manifolds are $S^3$-bundles over a quaternionic contact manifold equipped with a canonical conformal structure.
Moreover, the result from \cite{alt} also gives conformal holonomy  $G_{2,2}$ or  $Spin^+(4,3)$ for conformal structures of signature $(3,2)$ or $(3,3)$, respectively. Their geometric meaning has been revealed in \cite{hs1,nur,leinur}. In fact, these geometries can be equivalently characterized in terms of a generic 2-distribution on a 5-dimensional manifold, respectively a generic 3-distribution on a 6-dimensional manifold of split signature.\\
In Riemannian signature irreducibly acting conformal holonomy plays no role due to the fact that $O(1,n+1)$ admits no proper subgroup acting irreducibly on $\R^{p+1,q+1}$. Similarly, irreducibly acting subgroups of $Hol(M,c)$ can only occur in case $n$ even and in this case \cite{scal,leihoho} shows that they are exhausted by $SU(1, \frac{n}{2}) \subset SO(2,n)$.\\
Other initial results about conformal holonomy groups concerned the geometric meaning of a nontrivial subspace fixed by the standard action of $Hol(M,c)$ on $\R^{p+1,q+1}$. One finds that an invariant non-isotropic line corresponds to an Einstein metric in the conformal class off a singular set. In the isotropic case, a Ricci-flat metric occurs, cf. \cite{lst}. This mainly follows from fundamental properties of the covariant derivative $\nabla^{nc}$ on the conformal standard tractor bundle $\mathcal{T}(M)= \mathcal{P}^1 \times_{P} \R^{p+1,q+1}$ induced by $\omega^{nc}$ whose holonomy coincides with $Hol(M,c)$. Moreover, there is a conformal analogue of the de Rham decomposition theorem for metric holonomy groups. Concretely, a $Hol(M,c)$-invariant decomposition of $\R^{p+1,q+1}$ into nondegenerate factors of dimensions $r+1$ and $s+1$ corresponds for some open dense subset to a metric product of Einstein manifolds of dimensions $r$ and $s$ in the conformal class whose scalar curvatures satisfy a special linear relation, i.e. a special Einstein product in the sense of \cite{leihabil,arm,al2}. These results are the starting point for a complete classification of Riemannian holonomy groups which has been carried out in \cite{arm}.\\
\newline
In signatures higher than Riemannian, it is also possible for $Hol(M,c)$ to fix a totally lightlike subspace of dimension $>1$, which in terms of classification results turns out to be the most involved situation. Pseudo-Riemannian conformal holonomy groups fixing a totally lightlike subspace of dimension 2 have been studied in \cite{ln}: $Hol(M,c)$ fixes a totally lightlike nullplane if and only if on an open and dense subset of $M$, there is a metric $g \in c$ and a null line $L \subset TM$ such that $L$ is parallel wrt. $\nabla^g$ and $Ric^g(TM) \subset L$. The proof does not carry over to subspaces of dimension $>2$.
Consequently, we see that the most involved situation when dealing with non-irreducibly acting conformal holonomy occurs when $Hol_x(M,c)$ fixes a totally lightlike subspace of dimension $\geq 3$. Up to now there is no geometric description of this situation.

One aim of this article is to close this gap. Note that in case of non-irreducibly acting conformal holonomy with invariant subspace $V \subset \R^{p+1,q+1}$ one either has that $V$ is nondegenerate, which is covered by the previous review, or one can pass to a totally lightlike, nontrivial subspace $\widetilde{V} := V \cap V^{\bot}$ which is also fixed by the conformal holonomy representation. This case is solved locally in full generality here. It is in terms of tractors equivalent to the existence of a $\nabla^{nc}$-parallel and totally lightlike distribution in the standard tractor bundle $\mathcal{T}(M)$. In view of this, one has together with our main theorem a complete \textit{local} geometric description of conformal structures admitting non-irreducibly acting conformal holonomy. We prove:

\begin{satz} \label{gg}
If on a conformal manifold $(M,c)$ there exists a totally lightlike, $k$-dimensional parallel distribution $\mathcal{H} \subset \mathcal{T}(M)$, then there is an open and dense subset $\widetilde{M}$ of $M$ on which the totally lightlike distribution $L:=\text{pr}_{TM} (\mathcal{H} \cap \mathcal{I}^{\bot}_-) \subset TM$ canonically constructed out of $\mathcal{H}$ (as to be defined in section \ref{rt}) is of constant rank $k-1$ and integrable. Every point $x \in \widetilde{M}$ admits a neighbourhood $U=U_x \subset \widetilde{M}$ and a metric $g \in c_U$ such that on $U$:
\begin{equation}\label{3}
\begin{aligned}
&Ric^g(TU) \subset L, \\ 
&L \text{ is parallel wrt. }\nabla^g,\text{ i.e. }Hol_x(U,g)L_x \subset L_x.
\end{aligned}
\end{equation}
Conversely, let $(U,c)$ be a conformal manifold. Suppose that there is $g \in c$ and a $(k-1)$-dimensional totally lightlike distribution $L \subset TU$ such that (\ref{3}) holds. Then $L$ gives wrt. $g$ rise to a $k-$dimensional totally lightlike, parallel distribution $\mathcal{H}\stackrel{\Phi^g}{=} \begin{pmatrix} 0 \\ L \\ 0 \end{pmatrix} \oplus \text{span }\begin{pmatrix} 0 \\ 0 \\1 \end{pmatrix}$ in $\mathcal{T}(U)$, where $\mathcal{T}(U)$ is split wrt. $g$ via the map $\Phi^g$ as to be defined in section \ref{1rt}.
\end{satz}
Thus, one has a totally lightlike, parallel distribution in the standard tractor bundle if and only if one has locally a totally lightlike and parallel distribution of one dimension less in the tangent bundle with respect to some metric in the conformal class which additionally satisfies the curvature condition (\ref{3}). Up to now there is no complete classification of such metric holonomy groups .\\
Clearly, Theorem \ref{gg} generalizes the mentioned statements for an isotropic line or plane fixed by $Hol(M,c)$. Moreover, Theorem \ref{gg} also naturally generalizes results from \cite{nc} where the statement is proved under the additional condition that the totally lightlike distribution $\mathcal{H} \subset \mathcal{T}(M)$ arises from a decomposable, totally lightlike twistor $k-$form, by which we mean that there is a holonomy-invariant form $l_1 \wedge...\wedge l_k$ where the $l_i$ span a totally lightlike $k-$dimensional subspace in $\R^{p+1,q+1}$. Clearly, this space is then also holonomy-invariant. However, as elaborated on in \cite{ln}, in general not every holonomy-invariant totally lightlike $k-$dimensional subspace gives rise to a holonomy-invariant totally lightlike $k-$form. Thus we get the same geometric structures as discussed in \cite{nc} in the presence of totally lightlike twistor forms but under weaker assumptions.\\
\newline
The second aim of this article is to illustrate how Theorem \ref{gg} allows a classification of pseudo-Riemannian geometries admitting twistor spinors in certain low dimensions. For a space- and time oriented pseudo-Riemannian spin manifold $(M,g)$ with spinor bundle $S^g$, spinor covariant derivative $\nabla^{S^g}$ and Dirac operator $D^g$, they are given as solutions of the conformally-covariant twistor equation
\[ \nabla^{S^g}_X \ph + \frac{1}{n} X \cdot D^g \ph = 0 \text{ for } X \in TM. \]
Especially Riemannian and Lorentzian manifolds admitting twistor spinors have been well-studied and there many local geometric classification results, see \cite{bfkg,ha90,ha96,lei,bl,leihabil,bafe}. It has been observed in \cite{lei,baju,leihabil} that the twistor equation also admits a conformally invariant reinterpretation in terms of conformal Cartan geometries. In fact, there is a naturally associated vector bundle $\mathcal{S}$ for a conformal spin manifold $(M,c)$ of signature $(p,q)$ with fibre $\Delta_{p+1,q+1}$, the spinor module in signature $(p+1,q+1)$. On $\mathcal{S}$, a natural lift of the conformal Cartan connection $\omega^{nc}$ induces a covariant derivative such that parallel sections of $\mathcal{S}$ correspond to twistor spinors via a fixed metric $g \in c$. In other words, $(M,c)$ admits a twistor spinor for one - and hence for all - $g \in c$ iff the lift of $Hol(M,c)$ to the spin group $Spin^+(p+1,q+1)$ which double covers $SO(p+1,q+1)$ stabilizes a nonzero spinor. Using these Cartan techniques has lead to a complete local classification of Lorentzian conformal structures admitting twistor spinors in \cite{leihabil}.\\
A conformal holonomy group  in higher signature stabilizing a spinor is $G_{2,2}\subset SO^+(4,3)$ and leads to conformal structures of signature $(3,2)$ admitting twistor spinors. They have been intensively studied in \cite{hs1,nur,leinur,nur}, for instance. For these twistor spinors $\ph$, the distribution $H:= \text{ker }\ph = \{X \in TM \mid X \cdot \ph = 0 \}\subset TM$ is of constant rank 2 and turns out to be a generic 2-distribution, i.e. $[H,[H,H]]=TM$. Furthermore, the distribution ker $\ph$ associated to a twistor spinor $\ph$ in signature $(3,2)$ is generic iff $\langle \ph, D^g \ph \rangle \neq 0$. On the other hand, using the general machinery of parabolic geometries from \cite{cs}, \cite{hs1} shows that given any 5-dimensional manifold $M$ admitting an oriented, generic 2-distribution $H$, there is a canonical (Fefferman-type) construction of a conformal structure $[g]$ of signature $(3,2)$ on $M$ admitting a twistor spinor $\ph \in \Gamma(M,S^g)$ with $H= \text{ker }\ph$ and $Hol(M,c) \subset G_{2,2}$. A similar construction can be obtained in signature $(3,3)$: Twistor spinors satisfying $\langle \ph, D^g \ph \rangle{S^g} \neq 0$ are equivalent characterized in terms of generic 3-distributions on $TM$ (cf. \cite{br2}) and conformal holonomy reduction $Hol(M,c) \subset Spin^+(4,3) \subset SO^+(4,4)$.\\
Thus, to obtain a complete local geometric classification of twistor spinors in signatures $(3,2)$ and $(3,3)$, one has to consider those satisfying $\langle \ph, D^g \ph \rangle = 0$. In fact, this condition admits a nice reformulation in terms of conformally invariant tractor data on $\mathcal{S}$, and \cite{hs1} proves that the function $\langle \ph, D^g \ph \rangle$ is constant for any real twistor spinor in signature $(3,2)$ and $(3,3)$ and independent of $g \in c$. We then apply Theorem \ref{gg} to a natural distribution associated to the parallel spin tractor in $\mathcal{S}$ describing $\ph$ to prove:

\begin{satz} \label{tss}
Real twistor half-spinors $\ph$ in signature $(2,2)$ without zeroes and real twistor (half-)spinors without zeroes in signatures $(3,2)$ and $(3,3)$ satisfying that $\langle \ph , D^g \ph \rangle \equiv 0$ are locally conformally equivalent to parallel spinors (off a singular set). Their associated distributions ker $\ph:=\{X \in TM \mid X \cdot \ph = 0 \} \subset TM$ are integrable (off a singular set). The conformal holonomy representation in all these cases fixes a totally lightlike subspace of maximal dimension 3 resp. 4.
\end{satz}
As local normal forms for metrics admitting parallel real spinors in low-dimensional split-signatures are known from \cite{br,kath}, this statement classifies together with the results from \cite{hs1} all local conformal structures admitting real twistor spinors in the mentioned signatures.\\
\newline
This article is organized as follows: We recall how conformal structures can be described in terms of parabolic Cartan geometries in section \ref{1rt}. These preparations and notations enable us to prove Theorem \ref{gg} in section 3. We then focus on twistor spinors on conformal spin manifolds in section 4. Hereby, we first outline how they are equivalently characterized in terms of parallel spin tractors and then apply Theorem \ref{gg} to this setting to prove Theorem \ref{tss}.\\
\newline
\textbf{Acknowledgement} The author gladly acknowledges support from the DFG (SFB 647 - Space Time Matter at Humboldt University Berlin) and the DAAD (Deutscher Akademischer Austauschdienst / German Academic Exchange Service).


\section{Basic facts about conformal Cartan geometry} \label{1rt}
Let $G$ be an arbitrary Lie group with closed subgroup $P$. The $P-$bundle $G \rightarrow G/P$ together with the Maurer-Cartan form of $G$ serves as flat and homogeneous model for arbitrary Cartan geometries of type $(G,P)$. These are specified by the data $(\mathcal{G} \rightarrow M, \omega)$, where $M$ is a smooth manifold of dimension dim$(G/P)$, $\mathcal{G}$ is a $P-$principal bundle over $M$ and $\omega \in \Omega^1(\mathcal{G},\mathfrak{p})$, called the Cartan connection, is $Ad$-equivariant wrt. the $P-$action, reproduces the generators of fundamental vector fields and gives a pointwise linear isomorphism $T_u \mathcal{G} \cong \mathfrak{g}$. For detailed introduction to Cartan geometries, we refer to \cite{sharp,cs}. As a Cartan connection does not allow one to distinguish a connection in the sense of a right-invariant horizontal distribution in $\mathcal{G}$, it is convenient to pass to the enlarged principal $G-$bundle $\overline{\mathcal{G}}:= \mathcal{G} \times_P G$ on which $\omega$ induces a principal bundle connection $\overline{\omega}$, uniquely determined by $\iota^* \overline{\omega} = \omega$, where $\iota: \mathcal{G} \hookrightarrow \overline{\mathcal{G}}$ is the canonical inclusion. We then set $Hol_u(\mathcal{G},\omega):=Hol_{[u,e]}(\overline{\mathcal{G}},\overline{\omega}) \subset G$.\\

It is well-known that every conformal structure of signature $(p,q)$ with $n=p+q \geq 3$, i.e. an equivalence class $c$ of metrics differing by multiplication by a positive function, on a smooth manifold $M^n$ is equivalently encoded in a Cartan geometry $(\mathcal{P}^1 \rightarrow M, \omega^{nc})$ naturally associated to it via a construction called the first Prolongation of a conformal structure, cf. \cite{cs,baju,sharp,feh}. In this case, the group $G$ is given by $G=O(p+1,q+1)$ and the parabolic subgroup $P=Stab_{\R^+e_-}G$ is realized as the stabilizer of a lightlike ray $\R^+e_-$ under the natural $G-$action on $\R^{p+1,q+1}$. The homogeneous model is then given by the metric product $S^p \times S^q$ equipped with the conformal structure $[-g_p+g_q]$, where $g_p$ and $g_q$ are the round standard metric of the two spheres. One can also think of the homogeneous model as being a double cover of the projectivized lightcone in $\R^{p+1,q+1}$ equipped with a natural conformal structure on which $O(p+1,q+1)$ acts by conformal transformations. $\omega^{nc} \in \Omega^1(\mathcal{P}^1,\mathfrak{g})$ is called the normal conformal Cartan connection, and given $\mathcal{P}^1$, it is uniquely determined by the normalization condition $\partial^* \Omega^{nc}=0$, where $\Omega^{nc}:\mathcal{P}^1 \rightarrow Hom(\Lambda^2 \R^n, \mathfrak{so}(p+1,q+1))$ denotes the curvature of $\omega$ and $\partial^*$ denotes the Kostant codifferential, cf. \cite{cs}. \\
\newline
Given the standard action of $O(p+1,q+1)$ on $\R^{p+1,q+1}$, we obtain the associated standard tractor bundle $\mathcal{T}(M):=\mathcal{P}^1 \times_P \R^{p+1,q+1} = \overline{\mathcal{P}}^1 \times_G \R^{p+1,q+1}$ on which $\overline{\omega}^{nc}$ induces a covariant derivative $\nabla^{nc}$ which is metric wrt. the bundle metric $\langle \cdot, \cdot \rangle_{\mathcal{T}}$ on $\mathcal{T}(M)$ induced by the pseudo-Euclidean inner product on $\R^{p+1,q+1}$. We view $\nabla^{nc}$ as the conformal analogue of the Levi-Civita connection and define the conformal holonomy of $(M,c)$ for $x \in M$ to be
\[ Hol_x(M,c):=Hol_x(\nabla^{nc}) \subset O(\mathcal{T}_x(M), \langle \cdot , \cdot \rangle_{\mathcal{T}}) \cong O(p+1,q+1). \]
Obviously, $Hol(M,c) \cong Hol(\omega^{nc})$ as conjugated subgroups of $O(p+1,q+1)$. \\
\newline
By means of a metric in the conformal class, the conformally invariant objects introduced so far admit a more concrete description:
Any fixed $g \in c$ induces a so-called Weyl-structure in the sense of \cite{cs} and leads to a $O(p,q) \hookrightarrow O(p+1,q+1)$-reduction $\sigma^g: \mathcal{P}^g \rightarrow \mathcal{P}^1$. Here, $\mathcal{P}^g$ denotes the orthonormal frame bundle for $(M,g)$. Hereby, for the embedding $O(p,q) \hookrightarrow O(p+1,q+1)$ we split $\R^{p+1,q+1} = \R e_- \oplus \R^{p,q} \oplus \R e_+$, where $e_+$ is a lightlike vector such that $\langle e_-, e_+ \rangle = 1$ and the above sum is a direct sum of $O(p,q)-$modules. It follows that there is a $g-$metric splitting of the tractor bundle
\begin{align} 
\mathcal{T}(M) \stackrel{\Phi^g}{\cong} \underline{\R} \oplus TM \oplus \underline{\R} =: \mathcal{I}_- \oplus TM \oplus \mathcal{I}_+ \label{df}, 
\end{align}
under which tractors correspond to elements $(\alpha,X,\beta)$ and the tractor metric takes the form
 \begin{align} \langle (\alpha_1, Y_1, \beta_1), (\alpha_2, Y_2, \beta_2) \rangle_{\mathcal{T}} = \alpha_1 \beta_2 + \alpha_2 \beta_1 + g(Y_1,Y_2). \label{bum} \end{align}
The metric description of the tractor connection $\nabla^{nc}$, i.e. $\Phi^g \circ \nabla^{nc} \circ (\Phi^g)^{-1}$ is (cf. \cite{baju})
 \begin{align} \nabla_X^{nc} \begin{pmatrix} \alpha \\ Y \\ \beta \end{pmatrix} = \begin{pmatrix} X(\alpha) + K^g(X,Y) \\ \nabla_X^g Y + \alpha X - \beta K^g(X)^{\sharp} \\ X(\beta) - g(X,Y) \end{pmatrix}, \label{trad} \end{align}
where $K^g := \frac{1}{n-2} \cdot \left( \frac{scal^g}{2(n-1)}  \cdot g - Ric^g  \right)$ is the Schouten tensor.
Under a conformal change $\widetilde{g}=e^{2 \sigma}g$, the metric representation of tractors changes according to the map $\Phi^{\widetilde{g}} \circ (\Phi^{g})^{-1}$, given by (cf. \cite{baju})
\begin{align}
\begin{pmatrix} \alpha \\ Y \\ \beta \end{pmatrix} \mapsto \begin{pmatrix} \widetilde{\alpha} \\ \widetilde{Y} \\ \widetilde{\beta} \end{pmatrix}= \begin{pmatrix} e^{- \sigma} (\alpha - Y(\sigma) - \frac{1}{2}\beta ||\text{grad}^g \sigma ||^2_g \\ e^{- \sigma} (Y + \beta \text{grad}^g \sigma) \\ e^{\sigma} \beta \end{pmatrix}. \label{tra}
\end{align}

\section{Proof of Theorem \ref{gg}} \label{rt}
Theorem \ref{gg} later turns out to be a direct consequence of the following statement formulated on the level of tractors only:
\begin{Proposition} \label{ct}
Let $(M,c)$ be a conformal manifold of dimension $n\geq 3$ and let $\mathcal{H} \subset \mathcal{T}(M)$ be a totally lightlike distribution of dimension $k \geq 1$ which is parallel wrt. the Cartan connection $\nabla^{nc}$. Then there is an open, dense subset $\widetilde{M} \subset M$ such that for every point $x \in \widetilde{M}$ there is an open neighbourhood $U_x \subset \widetilde{M}$ and a metric $g \in c_{|U_x}$ such that wrt. the metric identification $\Phi^g$(cf. (\ref{df})) $\mathcal{H}$ is locally given by
\begin{align*}
\mathcal{H}_{|U_x} \stackrel{g}{=} \text{span }\left(\begin{pmatrix} 0 \\ K_1\\ 0 \end{pmatrix},...,\begin{pmatrix} 0 \\ K_{k-1}\\ 0 \end{pmatrix}, \begin{pmatrix} 0 \\ 0\\ 1 \end{pmatrix}\right)
\end{align*}
for lightlike vector fields $K_i \in \mathfrak{X}(U_x)$ which define a conformally invariant distribution \\
$L=\text{span }(K_1,...,K_{k-1})\subset TU_x$ of rank $k-1$ on $U_x$.
\end{Proposition}

\textbf{Proof. }If $k=1,2$, this statement is proved in \cite{baju} and \cite{ln}, respectively. Parts of the notations in this proof follow \cite{ln}. Consequently, we may assume that $k \geq 3$. As a preparation, consider for arbitrary $g \in c$ the map $\Phi^g: \mathcal{T}(M) \rightarrow \mathcal{I}_- \oplus TM \oplus \mathcal{I}_+=:\mathcal{T}(M)_g$. We set $\mathcal{I}_-=:\mathcal{I}$ and observe from the transformation formula (\ref{tra}) that this tractor null line which defines the conformal structure does not depend on the choice of $g \in c$. In this proof, we will for fixed $g$ always identify $\mathcal{T}(M)$ with $\mathcal{T}(M)_g$ without writing $\Phi^g$ explicitly. Moreover, we introduce the $g-$dependent projection $\text{pr}_{TM}:\mathcal{T}(M)\stackrel{g}{=}\mathcal{I}_- \oplus TM \oplus \mathcal{I}_+ \rightarrow TM$. Note however, that by (\ref{tra}), for every subbundle $\mathcal{V} \subset \mathcal{I}_-^{\bot} \stackrel{g}{=} \mathcal{I}_- \oplus TM$, the image $\text{pr}_{TM}\left(\mathcal{V} \right) \subset TM$ does not depend on the choice of $g \in c$. \\
We set $\mathcal{L}:=\mathcal{I}^{\bot} \cap \mathcal{H}$, where $\bot$ is taken wrt. the standard tractor metric.
By (\ref{bum}) we have that with respect to $g \in c$ it holds that $\mathcal{L} = \left\{ X \in \mathcal{H} \mid X = \begin{pmatrix} \alpha \\Y \\ 0 \end{pmatrix} \right\}$. It follows that $L:= \text{pr}_{TM} \mathcal{L} \subset TM$ is conformally invariant. With these introductory remarks in mind, the proof goes as follows:\\
\newline

\textbf{Step 1: }
We claim that there is an open, dense subset\footnote{In this proof, in order to keep notation short, whenever we restrict our considerations to an open, dense subset of $M$ we again call it $M$.} $\widetilde{M} \subset M$ such that $\text{rk }\mathcal{L}_{|\widetilde{M}} = k-1$:
Note that $\mathcal{L} \neq \{0 \}$ as otherwise $\mathcal{H}$ would have rank 1. Assume that there is an open set $U$ in $M$ on which $\text{rk }\mathcal{L}_{|U} = k$. We fix an arbitrary metric $g \in c$. By definition of $\mathcal{L}$, we have that $\mathcal{H} \cap \mathcal{I}^{\bot}=\mathcal{L} = \mathcal{H}$ on $U$ from which $\mathcal{H}_{|U} \subset \mathcal{I}^{\bot}_{|U} \stackrel{g}{=} (\mathcal{I}_- \oplus TM)_{|U}$ follows. Now let $\underline{L} \stackrel{g}{=} \begin{pmatrix} \rho \\ Y \\ 0 \end{pmatrix} \in \Gamma(\mathcal{H}_{U})$ be an arbitrary section of $\mathcal{H}$. As $\mathcal{H}$ is parallel, we must have that $\nabla^{nc}_X \underline{L} \in \Gamma(\mathcal{H}_{|U}) \subset \Gamma(\mathcal{I}^{\bot}_{|U})$ for all $X \in TU$. However, by (\ref{trad}) we get that 
\[ \nabla_X^{nc} \underline{L} = \begin{pmatrix} * \\ * \\ -g(X,Y) \end{pmatrix} \text{ }\forall X \in TU,\]
which means that $Y=0$ and $k=$rk $\mathcal{H}=1$.
Consequently, there is an open, dense subset (which we again call $M$) over which $0 < \text{rk } \mathcal{L} < k$. Now let $x \in M$ and fix a basis $L_1,...,L_s$ of $\mathcal{L}_x$, where $s=s(x) \leq k-1$. We may add tractors $Z_l = \begin{pmatrix} a_l \\ Y_l \\ 1 \end{pmatrix} \in \mathcal{T}_x(M)$ for $1 \leq l \leq k-s$ such that $L_1,...,L_s,Z_1,...,Z_{k-s}$ is a basis of of $\mathcal{H}_x$. We know that $k-s \geq 1$. If $k-s >1$ we may form new basis vectors $Z_1+Z_2$ and $Z_1-Z_2$. However, $Z_1-Z_2 \in \mathcal{L}_x$. Thus, $k-s=1$, which shows that rk $\mathcal{L}_x = k-1$.\\
\newline
\textbf{Step 2: }
We claim that also $L:=\text{pr}_{TM}\mathcal{L}$ has rank $k-1$ locally around each point $x \in M$. To this end, let $g \in c$ be arbitrary. Then we choose generators of $\mathcal{L}$ around $x$ such that locally $\mathcal{L} \stackrel{g}{=} \text{span } \left( \begin{pmatrix} \sigma_1 \\ \widetilde{K}_1 \\ 0 \end{pmatrix},...,\begin{pmatrix} \sigma_{k-1} \\ \widetilde{K}_{k-1} \\ 0 \end{pmatrix}\right)$. As $k>2$, we may assume that $\widetilde{K}_1(x) \neq 0$. We may then at the same time also assume that $\sigma_1(x) \neq 0$. Otherwise, we find $f \in C^{\infty}(M)$ with $\widetilde{K}_1(f)(x) \neq 0$ and consider the metric $\widetilde{g}=e^{2 f}g$ instead (cf. (\ref{tra})). Moreover, we may by multiplying the generators with nonzero functions assume that there is a neighbourhood $U$ of $x$ on which $\sigma_1 \equiv 1$ and $|\sigma_i| < 1$ for $i=2,...,k-1$. By applying elementary linear algebra to the generators, we then see that there are lightlike vector fields $K_i \in \mathfrak{X}(U)$ for $i=1,...,k-1$ with $K_1(x) \neq 0$ such that wrt. ${g}$ on $U$ 
\begin{align}
\mathcal{L} \stackrel{{g}}{=} \text{span } \left( \begin{pmatrix} 1 \\ {K}_1 \\ 0 \end{pmatrix}, \begin{pmatrix} 0 \\ {K}_2 \\ 0 \end{pmatrix},...,\begin{pmatrix} 0 \\ {K}_{k-1} \\ 0 \end{pmatrix}\right). \label{dd}
\end{align}
If $K_1$ was contained in the span of the $K_{i>1}$, we would have that $\begin{pmatrix} 1 \\ 0 \\ 0 \end{pmatrix} \in \mathcal{L} \subset \mathcal{H}$. However, as by Step 1 $\mathcal{H}$ must also contain a tractor of the form $\begin{pmatrix} a \\ X \\ 1 \end{pmatrix}$ not lying in $\mathcal{L}$, this contradicts $\mathcal{H}$ being totally lightlike. Consequently, there is an open neighbourhood of $x$ in $M$ such that the so constructed vectors $K_1,...,K_{k-1}$ are linearly independent and as pointwise
$L= \text{span} (K_1,...,K_{k-1})$ this shows that there is an open and dense subset of $M$ on which the rank of $L$ is maximal.\\
\newline
\textbf{Step 3: }
It follows directly from the various definitions that 
\begin{align}
\text{pr}_{TM} \left( \mathcal{L}^{\bot} \cap \mathcal{I}^{\bot} \right) = L^{\bot}. \label{ob1}
\end{align}
Moreover, note that $\mathcal{I} \subset \mathcal{L}^{\bot}$. By definition, $\mathcal{L} \subset \mathcal{H}$, from which $\mathcal{H}^{\bot} \subset \mathcal{L}^{\bot}$ follows. As by Step 1 $\mathcal{L}= \mathcal{H} \cap \mathcal{I}^{\bot}$ has codimension 1 in $\mathcal{H}$,  the line $\mathcal{I}$ cannot lie in $\mathcal{H}^{\bot}$, i.e. $\mathcal{H}^{\bot} \cap \mathcal{I} = \{0 \}$. A dimension count thus shows that
\begin{align}
\mathcal{L}^{\bot}=\mathcal{H}^{\bot} \oplus \mathcal{I}. \label{ob2}
\end{align}
(\ref{ob1}) and (\ref{ob2}) imply that
\[L^{\bot} = \text{pr}_{TM} \left( \mathcal{H}^{\bot} \cap \mathcal{I}^{\bot} \right).\]

\textbf{Step 4: }
In the setting of Step 2 we again fix $x \in M$, consider the local representation (\ref{dd}) of $\mathcal{L}$ wrt. some fixed $g$ around $x$ and set $L':=\text{span }(K_2,...,K_{k-1})$\footnote{Note that in contrast to $L$, the distribution $L'$ depends on the choice of $g \in c$.}. Both $L$ and $L'$ are integrable distributions. To see this, let $i,j \in \{2,...,k-1\}$. As $\mathcal{H}$ is parallel and totally lightlike we have that 
$\nabla^{nc}_{K_i} \begin{pmatrix} 0 \\ K_j \\ 0 \end{pmatrix} = \begin{pmatrix} K^g (K_i,K_j) \\ \nabla^g_{K_i}K_j \\ -g(K_i,K_j) \end{pmatrix} \in \Gamma(\mathcal{L})$. 
Switching the roles of $i$ and $j$ and taking the difference yields $\begin{pmatrix} 0 \\ \left[K_i,K_j\right] \\ 0 \end{pmatrix} \in \Gamma(\mathcal{L})$. Thus $\left[K_i,K_j\right] \in L'$. Similarly, one shows with the same argument that even
\begin{align}
[K_1,L'] \subset L'. \label{ic}
\end{align}
In particular, $L$ is integrable, too.\\
\newline
\textbf{Step 5: }
We now apply Frobenius Theorem to $L \subset TM$: For every (fixed) point $y$ of (an open and dense subset of ) $M$ we find a local chart $(U,\ph=(x_1,...,x_n))$ centered at $y$ with $\ph(U) = \{(x_1,...,x_n) \in \R^n \mid |x_i| < \epsilon \}$ such that the leaves $A_{c_k,...,c_n}=\{a \in U \mid x_k(a) = c_1,...,x_n(a)=c_n \} \subset U$ are integral manifolds for $L$ for every choice of $c_j$ with $|c_j| <\epsilon$. It holds that $L_U = \text{span }\left(\frac{\partial}{\partial x_1},...,\frac{\partial}{\partial x_{k-1}}\right)$ and moreover the coordinates may be chosen such that $K_1 = \frac{\partial}{\partial x_1}$ over $U$ (cf. \cite{war}). After applying some linear algebra to the generators of $L'$, where $L'$ is chosen wrt. some  $g \in c$ as in Step $4$ and restricting $U$ if necessary, we may assume that generators of $L'$ are given on $U$ by
\begin{align}
K_{i \geq 2} =  \alpha_i \frac{\partial}{\partial x_1} + \frac{\partial}{\partial x_i} \label{ki}
\end{align}
for certain smooth functions $\alpha_i \in C^{\infty}(U) \text{ for }i=2,...,k-1$. The integrability condition (\ref{ic}) implies that $\left[\frac{\partial}{\partial x_1},\alpha_i\frac{\partial}{\partial x_1} + \frac{\partial}{\partial x_i} \right]=\frac{\partial \alpha_i}{\partial x_1} \cdot \frac{\partial}{\partial x_1} \in L'$, giving that
\begin{align}
\frac{\partial \alpha_i}{\partial x_1}  = 0 \text{ for }i=2,...,k-1. \label{tt}
\end{align}
The integrability of $L'$ and (\ref{tt}) then yield that for $i,j=2,...,k-1$ 
\begin{align*}
[K_i,K_j]\stackrel{(\ref{ki}),(\ref{tt})}{=} \left( \frac{\partial \alpha_j}{\partial x_i}  - \frac{\partial \alpha_i}{\partial x_j} \right) \cdot \frac{\partial}{\partial x_1} \in L',
\end{align*}
from which by (\ref{ki}) follows that
\begin{align}
\frac{\partial \alpha_j}{\partial x_i}  - \frac{\partial \alpha_i}{\partial x_j} = 0 \text{ for }i,j=2,...,k-1. \label{jj}
\end{align}
For fixed $c_k,...,c_n$ as above we consider the submanifold $A_{c_k,...,c_n}$ and the differential form
\begin{align}
\alpha_{c_k,...,c_n} := - \sum_{i=1}^{k-1} \alpha_i dx_i \in \Omega^1\left(A_{c_k,...,c_n} \right),
\end{align}
where the $\alpha_{i \geq 2}$ are restrictions of the functions appearing in (\ref{ki}) to $A_{c_k,...,c_n}$ and we set $\alpha_1 \equiv -1$. (\ref{tt}) and (\ref{jj}) precisely yield that $d\alpha_{c_k,...,c_n} = 0$. Whence, there exists by the Poincar{\'e} Lemma (applied to a sufficiently small simply-connected neighbourhood) a unique $\sigma_{c_k,...,c_n} \in C^{\infty}\left(A_{c_k,...,c_n} \right)$ with $\sigma_{c_k,...,c_n}(\ph^{-1}(0,...,0,c_k,...,c_n)) = 0$ and $\alpha_{c_k,...,c_n} = d\sigma_{c_k,...,c_n}$, which translates into
\begin{align*}
\frac{\partial \sigma_{c_k,...,c_n}}{\partial x_1}  &= 1, \\
\frac{\partial  \sigma_{c_k,...,c_n}}{\partial x_i} &= - \alpha_i \text{ for }i=2,...,k-1.
\end{align*}
We define $\sigma \in C^{\infty}(U)$ via $\sigma(\ph^{-1}(x_1,....,x_n)):=\sigma_{x_k,...,x_n}(\ph^{-1}(x_1,...,x_{n}))$ and observe that on $U$
\begin{equation}\label{hhh}
\begin{aligned}
\frac{\partial \sigma}{\partial x_1} &= 1, \\
\frac{\partial \sigma}{\partial x_i}  &= - \alpha_i \text{ for }i=2,...,k-1. 
\end{aligned}
\end{equation}
\textbf{Step 6: }
The construction of the generators $K_i$ (\ref{ki}) and the properties (\ref{hhh}) of $\sigma$  imply that on $U$ we have $K_1(\sigma)=1$ and $K_i(\sigma) = 0$ for $i=2,...,k-1$. We now consider the rescaled metric ${\widetilde{g}}=e^{2\sigma}{g}$ on $U$. The transformation formula (\ref{tra}) and (\ref{hhh}) then show that wrt. this metric $\mathcal{L}$ is given by
\begin{align}
\mathcal{L}_U = \text{span} \left(\begin{pmatrix} 0 \\ K_1 \\ 0 \end{pmatrix},...,\begin{pmatrix} 0 \\ K_{k-1} \\ 0 \end{pmatrix} \right). \label{thepro}
\end{align}
\textbf{Step 7: }
Let $g \in c$ be any local metric on $U \subset M$ for which (\ref{thepro}) holds. We may add one generator $\begin{pmatrix} \beta \\ K \\1 \end{pmatrix} \in \Gamma(U,\mathcal{H})$ such that pointwise (wrt. ${g}$) $\mathcal{H}= \mathcal{L} \oplus \text{span }\begin{pmatrix} \beta \\ K \\1 \end{pmatrix}$. It follows that $K \in L^{\bot}$ as $\mathcal{H}$ is totally lightlike. By step $3$ there exists a smooth function $b$ on $U$ with $K= \text{pr}_{TM} \begin{pmatrix} b \\ K \\0 \end{pmatrix}$ and $\begin{pmatrix} b \\ K \\0 \end{pmatrix} \in \mathcal{H}^{\bot}$. As $\mathcal{H}$ is lightlike, (\ref{bum}) yields that $\beta + {g}(K,K)=0$ as well as $b+{g}(K,K)=0$, i.e. $b=\beta$. Therefore we have that $\begin{pmatrix} 0 \\ 0 \\1 \end{pmatrix} \in \mathcal{H}^{\bot}$ over $U$. However, this implies that $b=\beta = 0$ and we obtain 
\begin{align}
\mathcal{H}_U \stackrel{g}{=} \text{span} \left(\begin{pmatrix} 0 \\ K_1 \\ 0 \end{pmatrix},...,\begin{pmatrix} 0 \\ K_{k-1} \\ 0 \end{pmatrix}, \begin{pmatrix} 0 \\ K \\ 1 \end{pmatrix} \right). \label{rer}
\end{align}
In the following steps we will improve the fixed metric $g$ satisfying (\ref{thepro}) within the conformal class further in such a way that $K$ can be chosen to be zero. This goes as follows: Let $X \in \Gamma(L)$ be an arbitrary, nonzero section. We have for $Y \in TM$ that
\begin{align*}
\nabla^{nc}_Y \begin{pmatrix} 0 \\ X \\ 0 \end{pmatrix} \stackrel{(\ref{trad})}{=} \begin{pmatrix} * \\ \nabla_Y^g X \\ -g(X,Y) \end{pmatrix} \in \Gamma(\mathcal{H}),
\end{align*} 
yielding $\nabla^g_T X \in \Gamma(L)$ for $T \in X^{\bot}$ and for perpendicular directions
\begin{align}
\nabla^g_Z X = l-K \label{gy}
\end{align}
 for some $l \in L$, where $g(X,Z)=1$. Thus, if $g \in c$ can be chosen such that (\ref{thepro}) and additionally $\nabla_Y^g X \in \Gamma(L)$ hold for every $Y \in TM$, it holds that $K \in L$ and we can obviously rearrange the generators in (\ref{rer}) such that the Proposition follows.
Steps 8 and 9 are a preparation for the construction of this desired metric.
\newline
\textbf{Step 8: }
Wrt. $g \in c$ a metric satisfying (\ref{thepro}), let $\begin{pmatrix} \rho \\ V \\ 0 \end{pmatrix} \in \Gamma \left(\mathcal{H}^{\bot} \cap \mathcal{I}^{\bot}\right)$. Further, let $X \in \Gamma(L)$ be nonzero and let $Z$ be a vector field with $g(X,Z)=1$. As $\mathcal{H}$ is parallel and lightlike, we have 
\begin{align}
0 = \langle \nabla_Z^{nc} \begin{pmatrix}0 \\ X \\ 0 \end{pmatrix}, \begin{pmatrix} \rho \\ V \\ 0 \end{pmatrix} \rangle_{\mathcal{T}}= -\rho +g(\nabla_Z^g X,V). \label{gar1}
\end{align}
Let $U \in L^{\bot}$. Further differentiation yields
\begin{align*}
\nabla_U^{nc} \nabla_Z^{nc} \begin{pmatrix}0 \\ X \\ 0 \end{pmatrix} = \begin{pmatrix} * \\ \nabla_U^g \nabla_Z^g X + K^g(X,Z) \cdot U +K^g(U)^{\sharp}  \\ -g(U,\nabla^g_Z X)\end{pmatrix} \in \Gamma(\mathcal{H}).
\end{align*}
Pairing with $\begin{pmatrix} \rho \\ V \\ 0 \end{pmatrix} \in \Gamma \left(\mathcal{H}^{\bot} \cap \mathcal{I}^{\bot}\right)$ leads to
\begin{align}
0 = -\rho \cdot g(U,\nabla^g_Z X) + g(\nabla_U^g \nabla_Z^g X,V)+K^g(U,V)+K^g(X,Z)\cdot g(U,V). \label{gar2}
\end{align}
It follows from (\ref{gar1}) and (\ref{gar2}) that the bilinear form
\begin{align}
\Gamma(L^{\bot}) \times \Gamma(L^{\bot}) \ni (U,V) \mapsto g(\nabla_U^g \nabla_Z^g X,V) \label{syga}
\end{align}
is symmetric.\\
\newline
\textbf{Step 9: }
Let $g \in c$ be a metric satisfying (\ref{thepro}). Locally, we have that $L=\text{span} \left(K_1,...,K_{k-1}\right)$ and $L^{\bot}=\text{span} \left(K_1,...,K_{k-1},E_1,...,E_l\right)$, where $l=n-2k+2$ and the $E_i$ are vector fields on $U \subset M$ which are orthogonal to the $K_i$ and satisfy $g(E_i,E_j)= \pm 1 \cdot \delta_{ij}$. $L^{\bot}$ is an integrable distribution: By Step 3 there exist functions $\rho_j$ such that $E_j$ is the projection of $\begin{pmatrix} \rho_j \\ E_j \\ 0 \end{pmatrix} \in \Gamma \left(\mathcal{H}^{\bot} \cap \mathcal{I}^{\bot}\right)$ to $TM$ for $j=1,...,n-2k+2$. As also $\mathcal{H}^{\bot}$ is parallel, (\ref{trad}) yields for $i \neq j$ that
\begin{align*}
\nabla_{E_i}^{nc} \begin{pmatrix} \rho_j \\ E_j \\ 0 \end{pmatrix} = \begin{pmatrix} * \\ \nabla^g_{E_i}E_j +\rho_j E_i \\ 0 \end{pmatrix} \in \Gamma \left(\mathcal{H}^{\bot} \cap \mathcal{I}^{\bot}\right).
\end{align*} 
It follows that $\nabla^g_{E_i}E_j \in L^{\bot}$. With the same argumentation, one finds that also $\nabla^g_{E_i} K_j, \nabla^g_{K_j} E_i \in L^{\bot}$ for $i=1,...,n-k+1$, $j=1,...k-1$. This yields together with integrability of $L$ and torsion-freeness of $\nabla^g$ the integrability of $L^{\bot}$.\\
\newline
\textbf{Step 10: }
As $L \subset L^{\bot}$ and both are integrable distributions, we can by Frobenius Theorem (cf. Step 5) applied first to $L^{\bot}$ and then to each leaf of $L^{\bot}$ find around every point local coordinates \[(U,\ph=(x_1,...,x_{k-1},y_1,...,y_{n-2k+2},z_1,...,z_{k-1}))\] 
such that $(x_1,...,x_{k-1})$ parametrizes integral manifolds for $L$ and $(x_1,...,x_{k-1},y_1,...,y_{n-2k+2})$ parametrizes integral manifolds for $L^{\bot}$. \\
Let $\sigma \in C^{\infty}(U)$ be an arbitrary function depending on $(y_1,...,y_{n-2k+2},z_1,...,z_{k-1})$ only and set $\widetilde{g}=e^{2 \sigma} g$. Clearly, $X(\sigma)=0$ for every $X \in L=\text{span}\left(\frac{\partial}{\partial x_1},...,\frac{\partial}{\partial x_{k-1}}\right)$ and thus the tractor $\begin{pmatrix} 0 \\ X \\ 0 \end{pmatrix} \in \Gamma(\mathcal{L})$ is wrt. $\widetilde{g}$ given by (cf.(\ref{tra})) $\begin{pmatrix} 0 \\ \widetilde{X} \\ 0 \end{pmatrix} \in \Gamma(\mathcal{L})$ for some $\widetilde{X} \in L$. This means that also $\widetilde{g}$ satisfies (\ref{thepro}).\\
We set $X := \frac{\partial}{\partial x_1} \in \Gamma(L)$ and fix a vector field $Z$ such that $g(X,Z)=1, g( \frac{\partial}{\partial x_{i>1}},Z)=g( \frac{\partial}{\partial y_{j}},Z)=0$. We want to show that 
\begin{align}
g(\nabla^{\widetilde{g}}_Z X,Y) = 0 \text{ for every }Y \in L^{\bot}, \label{wantsh}
\end{align}
from which $\nabla^{\widetilde{g}}_Z X \in \Gamma(L)$ follows. To this end, we calculate with the well-known transformation formula $\nabla^{\widetilde{g}}_B A = \nabla_B^gA +d\sigma(B) A + d\sigma(A) B - g(A,B) \cdot  \text{grad}^g \sigma$ for the Levi-Civita connection that
\begin{equation} \label{depo2}
\begin{aligned}
g(\nabla^{\widetilde{g}}_Z X,Y)&=g(\nabla^g_Z X,Y)+d\sigma(Z)\underbrace{g(X,Y)}_{=0}+\underbrace{d\sigma(X)}_{=0}g(Z,Y)-g(\text{grad}^g\sigma,Y)\\
&=\left(g(\nabla^g_Z X,\cdot) - d\sigma \right)(Y),
\end{aligned}
\end{equation}
where $Y \in L^{\bot}$. On the other hand, we calculate for $U,V \in \Gamma(L^{\bot})$
\begin{equation} \label{depo}
\begin{aligned}
d(g(\nabla^g_Z X, \cdot))(U,V) &= U(g(\nabla_Z^gX,V))-V(g(\nabla^g_ZX,U))-g(\nabla^g_ZX,[U,V]) \\
&=g(\nabla_U^g\nabla^g_ZX,V)-g(\nabla_V^g\nabla^g_ZX,U) \stackrel{(\ref{syga})}{=} 0.
\end{aligned}
\end{equation}
To evaluate this further, we introduce $\theta:=g(\nabla^g_Z X, \cdot) \in \Omega^1(U)$. As moreover $g(\nabla^g_Z X, l)=0$ for every $l \in L$ (cf. (\ref{gy})), there exist local functions $\alpha_i, \beta_j \in C^{\infty}(U)$ such that $\theta = \sum_i \alpha_i dy_i + \sum_j \beta_j dz_j$. Let us define $\widetilde{\theta}:=\sum_i \alpha_i dy_i$ and let $\widetilde{\theta}_{A_{c_1,...,c_{k-1}}}$ denote its restriction to the leaf $A_{c_1,...,c_{k-1}}:=\{ \ph(x_1,...,x_{k-1},y_1,...,y_{n-2k+2},c_1,...,c_{k-1}) \mid c_i = \text{const}. \}$ of $L^{\bot}$. Obviously, (\ref{depo}) is equivalent to $d\left(\widetilde{\theta}_{A_{c_1,...,c_{k-1}}}\right)=0$ for all $c_i$. Thus, by applying the Poincar{\'e} Lemma again on a sufficiently small neighbourhood, we conclude that there are unique $\gamma_{c_1,...,c_{k-1}} \in C^{\infty}(A_{c_1,...,c_{k-1}})$ such that $\gamma_{c_1,...,c_{k-1}}(\ph^{-1}(0,...,0,c_1,...,c_{k-1}))=0$ and $d\gamma_{c_1,...,c_{k-1}}=\widetilde{\theta}_{A_{c_1,...,c_{k-1}}}$. We now specify $\sigma \in C^{\infty}(U)$ by setting 
\[ \sigma(\ph^{-1}(x_1,....,y_{n-2k+2},z_1,...,z_{k-1})):= \gamma_{z_1,...,z_{k-1}}(\ph^{-1}(x_1,....,y_{n-2k+2},z_1,...,z_{k-1})).
\]
This construction yields for $Y \in L^{\bot}$
\[ d\sigma(Y)=\widetilde{\theta}(Y)=\theta(Y)=g(\nabla^g_Z X, Y). \]
Letting $Y=\frac{\partial}{\partial x_i}$ and using $\nabla^g_Z X \in \Gamma(L^{\bot})$, cf. (\ref{gy}), yields $\frac{\partial \sigma}{\partial x_i}=0$, i.e. $\sigma$ does not depend on $(x_1,...,x_{k-1})$. Consequently, we get from (\ref{depo2}) for this choice of $\sigma$ that (\ref{wantsh}) holds.  However, as remarked at the end of Step 7, this already proves the Proposition.
$\hfill \Box$ \\

We study some consequences. In the setting of Proposition \ref{ct} we have that $\mathcal{H}$ is parallel iff $\mathcal{H}^{\bot}$ is parallel. Locally, we have wrt. the metric $g$ and the distribution $L$ appearing in Proposition \ref{ct} that $\mathcal{H}^{\bot}=\text{span} \left( \begin{pmatrix} 0 \\ X \\ \tau \end{pmatrix} \mid X \in L^{\bot} \right)$. It follows that $\mathcal{H}^{\bot}$ is parallel iff

\begin{align*}
\nabla^{nc}_Y \begin{pmatrix} 0 \\ X \\ \tau \end{pmatrix} = \begin{pmatrix} K^g(X,Y) \\ \nabla^g_Y X - \tau K^g(Y) \\ Y(\tau)-g(X,Y) \end{pmatrix} \in \Gamma(U,\mathcal{H}^{\bot})
\end{align*}
for all $X \in \Gamma(U,L^{\bot})$ and $Y \in \mathfrak{X}(U)$. Clearly, this is equivalent to parallelism of $L$ and $K^g(X,Y)=0$ for all $X \in L^{\bot}$, i.e. $K^g(TU) \subset L$. Together with the next Lemma, these two conditions are equivalent to parallelism of $L$ and $Ric^g(TU) \subset L$.

\begin{Lemma}
Assume that for a pseudo-Riemannian manifold $(M,g)$ one has a nontrivial totally lightlike $(k-1)-$dimensional distribution $L \subset TM$ for which $K^g(TM) \subset L$. Then $\text{scal}^g=0$.
\end{Lemma}

\textbf{Proof. }
For fixed $x \in M$ we introduce a basis $(X_1,...,X_{k-1},X'_1,...,X'_{k-1},E_1,...,E_l)$ of $T_xM$, where $L_x=\text{span}\{X_1,...,X_{k-1}\}$, $g(X_i,X'_j)=\delta_{ij}, g(X'_i,X'_j)=0$, $g(E_i,E_j)=\epsilon_i \delta_{ij}$ and $g(E_i,X_j^{(')})=0$. It follows that
\begin{align}
\text{scal}^g(x)=2 \sum_{j=1}^{k-1}Ric^g(X_j,X'_j) + \sum_{i=1}^l\epsilon_i Ric^g(E_i,E_i). \label{frodo}
\end{align}
By definition of the Schouten tensor, we have that $Ric^g = \frac{1}{2(n-1)} \cdot \text{scal}^g - (n-2) \cdot K^g$. Inserting this into (\ref{frodo}) yields
\begin{align*}
\text{scal}^g(x)=&\frac{1}{n-1}(k-1) \cdot \text{scal}^g(x) - 2(n-2) \cdot \sum_{j=1}^{k-1} \underbrace{K^g(X_j,X'_j)}_{=0} \\
&+ \frac{1}{2(n-1)} \cdot (n-2(k-1)) \cdot \text{scal}^g(x) -(n-2)\sum_{j=1}^l \epsilon_i \underbrace{K^g(E_i,E_i)}_{=0} \\
=&\frac{n}{2(n-1)}\cdot \text{scal}^g(x),
\end{align*}
i.e. $\text{scal}^g(x)=0$.
$\hfill \Box$ \\
\newline
Finally, we have seen in the proof of Proposition \ref{ct} that $L= \text{pr}_{TM} \mathcal{L} \subset TM$ is a well-defined distribution of constant rank on $\widetilde{M} \subset M$ open and dense. As $L$ is on $\widetilde{M}$ \textit{locally} parallel wrt. certain metrics in the conformal class, this implies by the torsion-freeness of $\nabla^g$ as a \textit{global} consequence that $L$ is integrable on $\widetilde{M}$. Thus, altogether we have proved Theorem \ref{gg}.

\begin{bemerkung}
It is a common feature of all statements about reducible conformal holonomy that one always has to leave out a certain set of singular points, i.e. restrict to some open and dense subset $\widetilde{M} \subset M$, as was also necessary in the proof of Proposition \ref{gg}. The deeper reason for this has recently been discovered in \cite{cgh}, and it is closely related to so called curved orbit decompositions of arbitrary Cartan geometries. At least in the case of a holonomy-invariant line, the reference can with this method also describe the geometry of the singular set $M \backslash \widetilde{M}$.
\end{bemerkung}

Let us study the local geometries occurring in Theorem \ref{gg} in more detail: Pseudo-Riemannian geometries admitting parallel, totally lightlike distributions are called Walker manifolds and have been studied in \cite{wal}, for instance. Let us call a pseudo-Riemannian manifold $(M,g)$ admitting a parallel, totally lightlike distribution $L \subset TM$  of rank $r$, satisfying additionally that $Ric^g(TM) \subset L$ a Ricci-isotropic pseudo-r-Walker manifold.\\
In general, for every $n-$dimensional Walker manifold $(M,g)$ with parallel, $r-$dimensional, totally lightlike distribution $L \subset TM$, there are locally around each point coordinates $(x_1,...,x_n)$ such that wrt. the basis $\left(\frac{\partial}{\partial x_1},...,\frac{\partial}{\partial x_n} \right)$ the metric tensor reads (cf. \cite{wal})
\begin{align*}
g_{ij} = \begin{pmatrix} 0 & 0 & Id_r \\ 0 & A & H \\ Id_r & H^T & B \end{pmatrix},
\end{align*}
where $A$ is a symmetric $(n-2r) \times (n-2r)$ matrix, $B$ is a symmetric $r \times r$ matrix and $H$ is a $(n-2r) \times r$ matrix. Moreover, $A$ and $H$ do not depend on $(x_1,...,x_r)$, and in these coordinates, $L$ is given by
\begin{align*}
L = \text{span}\left( \frac{\partial}{\partial x_1},...,\frac{\partial}{\partial x_r} \right).
\end{align*}

\begin{beispiel}Let $\ph$ be a parallel spinor on a pseudo-Riemannian spin manifold $(M,g)$. Then $L:=\text{ker }\ph= \{X \in TM \mid X \cdot \ph = 0 \} \subset TM$ is totally lightlike and parallel. $Ric^g(X) \cdot \ph = 0$ as known from \cite{ba81} translates into $Ric^g(TM) \subset L$. For small dimensions all Ricci-isotropic pseudo-$r$-Walker metrics arising in this way have been classified in \cite{br}. The orbit structure of $\De$ encodes which values for $r=$dim $L$ are possible in these cases.
\end{beispiel}

\section{Application to twistor spinors}

Let us briefly recall the twistor equation on spinors, starting with some spinor-algebraic facts (cf. \cite{ba81,har,lm}). By $SO^+(p,q) $ we denote the connected component of the orthogonal group in signature $(p,q)$, with $n=p+q \geq 3$, which is double-covered by $Spin^+(p,q)$, the connected component of the spin group, by means of a smooth group homomorphism $\lambda$. Let moreover $\Delta_{p,q}^{\R}$ denote the real spinor module in signature $(p,q)$, on which both $Spin^+(p,q)$ and $\R^{p,q}$ act by the standard spinor representation and Clifford multiplication, respectively, and it holds that $x \cdot x = -||x||^2$ for $x \in \R^{p,q}$. $\Delta_{p,q}^{\R}$ is equipped with a nondegenerate inner product $\langle \cdot, \cdot \rangle_{\Delta_{p,q}^{\R}}$ which is symmetric or symplectic depending on $q-p$ mod 4 and invariant under the identity component $Spin^+(p,q)$.\\
\newline
Turning to geometry, let $(M,g)$ be a space- and time oriented pseudo-Riemannian spin manifold of signature $(p,q)$, where $n=p+q \geq 3$. In this case the orthonormal frame bundle $\mathcal{P}^g$ can always be reduced further to the $SO^+(p,q)-$bundle of space-and time oriented pseudo-orthonormal frames, denoted by the same symbol. We fix a spin structure $(\mathcal{Q}^g, f^g : \mathcal{Q}^g \rightarrow \mathcal{P}^g)$, i.e. a $\lambda-$reduction of $\mathcal{P}^g$. The real spinor bundle is given by $S^g := \mathcal{Q}^g \times_{Spin^+(p,q)} \Delta_{p,q}^{\R}$ and it is equipped with Clifford multiplication by elements of $TM$ and an inner product $\langle \cdot, \cdot \rangle_{S^g}$. The Levi Civita connection on $\mathcal{P}^g$ lifts to a connection on $\mathcal{Q}^g$ via $f^g$, which in turn induces on $S^g$ a spinor covariant derivative, 
\[ \nabla^{S^g}: \Gamma(S^g) \rightarrow \Gamma(T^*M \otimes S^g). \]
Superposition of $\nabla^{S^g}$ with Clifford multiplication defines the Dirac operator $D^g:\Gamma(S^g) \rightarrow \Gamma(S^g)$, whereas superposition of $\nabla^{S^g}$ with projection onto the kernel of Clifford multiplication in $T^*M \otimes S^g$ defines a complementary operator $P^g$, called the Penrose- or twistor operator. Elements in its kernel are called twistor spinors, or conformal Killing spinors, and they are equivalently characterized as solutions of the conformally covariant twistor equation
\[ \nabla_X^{S^g} \ph + \frac{1}{n} X \cdot D^g \ph = 0 \text{ for }X \in TM. \]
Under a conformal change $\widetilde{g}=e^{2 \sigma}g$, there is a natural identification $\widetilde{}: S^g \rightarrow S^{\widetilde{g}}$ (cf. \cite{ba81,bfkg}), and it holds that $\ph \in \text{ker }P^g$ iff $e^{\frac{\sigma}{2}} \widetilde{\ph} \in \text{ker }P^{\widetilde{g}}$.
Conformal Cartan geometry allows a conformally invariant construction of $P^g$. Suppose that $(M,c)$ is space- and time oriented and spin for one - and hence for all - $g \in c$.  The construction from section 2 admits finer underlying structures: $(\mathcal{P}^1, \omega^{nc})$ can be reduced to a Cartan geometry of type $(G^+,P^+)$, denoted by the same symbol. It lifts to a conformal spin Cartan geometry $(\mathcal{Q}^1, \widetilde{\omega}^{nc})$ of type $(Spin^+(p+1,q+1),\widetilde{P}^+:=\lambda^{-1}(P^+))$ with associated spin tractor bundle
\[ \mathcal{S}:= \mathcal{Q}^1 \times_{\widetilde{P}^+} \Delta_{p+1,q+1}^{\R}, \]
on which $\mathcal{T}(M)$ acts by fibrewise Clifford multiplication and $\widetilde{\omega}^{nc}$ induces a covariant derivative $\nabla^{\mathcal{S}}$ on $\mathcal{S}$.
Fixing a metric $g \in c$ leads to a $Spin^+(p,q) \hookrightarrow Spin^+(p+1,q+1)$-reduction $\widetilde{\sigma}^g: \mathcal{Q}^g \rightarrow \mathcal{Q}^1$ which covers $\sigma^g$. We let $\overline{\mathcal{Q}}^1_+$ denote the enlarged $Spin^+(p+1,q+1)$-principal bundle and use $g$ to identify $\mathcal{S}(M) \cong Q^g_+ \times_{Spin^+(p,q)} \Delta^{\R}_{p+1,q+1}$.
However, $\Delta^{\R}_{p+1,q+1} = Ann(e_+) \oplus Ann(e_-)$ as $Spin^+(p,q)$-representations, where $Ann(e_{\pm}) = \{ w \in \Delta^{\R}_{p+1,q+1} \mid x \cdot w = 0 \}$ are two copies of $\Delta_{p,q}$, leading to natural $g-$dependent projections $\text{proj}_{\pm}^g: \mathcal{S}(M) \rightarrow {\mathcal{Q}}^g \times_{Spin^+(p,q)}Ann(e_{\pm})$ and the $g-$metric identification
\begin{align}
\widetilde{\Phi}^g: \mathcal{S}(M) \rightarrow S^g(M) \oplus S^g(M) \label{gdg}
\end{align}
One calculates that under (\ref{gdg}), $\nabla^{nc}$ is given by the expression (cf. \cite{baju})
\begin{align}
\nabla^{nc}_X \begin{pmatrix} \ph  \\ \phi \end{pmatrix} = \begin{pmatrix} \nabla_X^{S^g} & -X \cdot \\ \frac{1}{2}K^g(X) \cdot & \nabla^{S^g}_X \end{pmatrix} \begin{pmatrix} \ph  \\ \phi \end{pmatrix}. \label{h5}
\end{align}
As every twistor spinor $\ph \in \text{ker }P^g$ satisfies $\nabla^{S^g}_X \ph = \frac{n}{2}K(X) \cdot \ph$, cf. \cite{bfkg}, this yields a reinterpretation of twistor spinors in terms of conformal Cartan geometry. 
Namely for any $g \in c$, the vector spaces ker $P^g$ and parallel sections in $\mathcal{S}(M)$ wrt. $\nabla^{nc}$ are naturally isomorphic via
\begin{align}
 \text{ker }P^g \rightarrow \Gamma(S^g(M) \oplus S^g(M)) \stackrel{\left(\widetilde{\Phi}^g\right)^{-1}}{\cong } \Gamma(\mathcal{S}(M))\text{,   } 
 \ph \mapsto \begin{pmatrix} \ph \\  -\frac{1}{n}D^g \ph \end{pmatrix} \stackrel{\left(\widetilde{\Phi}^g\right)^{-1}}{\mapsto} \psi \in Par(\mathcal{S}_{\mathcal{T}}(M), \nabla^{nc}), \label{rei}
 \end{align}
 i.e. a spin tractor $\psi \in \Gamma(\mathcal{S}(M))$ is parallel iff for one - and hence for all -  $g \in c$ it holds that $\ph:= \widetilde{\Phi}^g(\text{proj}_+^g \psi) \in \text{ker }P^g$ and in this case $D^g \ph = -n \cdot \widetilde{\Phi}^g(\text{proj}_-^g \psi)$.\\

In combination with the previous results we can detect whether a twistor spinor $\ph$ is locally conformally equivalent to a parallel spinor on the level of Cartan geometries as follows:
Let $\psi \in \Gamma(\mathcal{S}(M))$ be a parallel spin tractor. We set \begin{align}\text{ker }{\psi(x)}:= \{ v \in \mathcal{T}_x(M) \mid  v \cdot \psi(x) = 0 \}.\label{ker} \end{align} Performing this for every point yields a totally lightlike distribution $\text{ker }{\psi} \subset \mathcal{T}(M)$. It is moreover parallel wrt. $\nabla^{nc}$, and henceforth of constant rank, since $Y \in \Gamma(\text{ker }\psi)$ and $X \in \mathfrak{X}(M)$ implies that $0=\nabla^{nc}_X (Y \cdot \psi) = \left( \nabla^{nc}_X Y \right) \cdot \psi$. 
Consequently, every parallel spin tractor $\psi$  naturally gives rise to a -possibly trivial- distinguished totally lightlike subspace fixed by the holonomy representation, $Hol_x(M,c) \text{ ker }\psi(x) \subset \text{ ker }\psi(x)$.
In complete analogy, if $\ph \in \Gamma(S^g)$ is parallel wrt. some $g \in c$, we get a totally lightlike, parallel distribution ker ${\ph} \subset TM$.

\begin{Proposition} \label{tms}
If $\psi \in \Gamma(\mathcal{S}(M))$ is a parallel spin tractor with ker ${\psi} \neq \{0 \}$, then there is an open and dense subset $\widetilde{M} \subset M$ such that on $\widetilde{M}$ the associated twistor spinor $\ph:=\widetilde{\Phi}^g \left( {proj}^g_+\psi \right)$ is locally conformally equivalent to a parallel spinor. Moreover, on $\widetilde{M}$ the distribution ker $\ph$ is of constant rank and integrable.
\end{Proposition}

\textbf{Proof. }Proposition \ref{ct} applied to $\mathcal{H}=$ker $\psi$ yields the desired $\widetilde{M}$ and for $x \in \widetilde{M}$ a neighbourhood $U$ and a local metric $g=g_{U}\in c_U$ such that wrt. $g$ we have $s_+=\begin{pmatrix} 0 \\ 0 \\ 1 \end{pmatrix} \in \text{ker }{\psi_{|U}}$. If we decompose $\psi$ on $U$ wrt. $g$ as in (\ref{rei}), i.e. $\psi_{|U} = \left[\left[\widetilde{\sigma}^g(\widetilde{u}),e \right], e_- \cdot w + e_+ \cdot w \right]$ for some function $w:U \rightarrow \Delta_{p+1,q+1}$ and a local section $\widetilde{u}:U \rightarrow \mathcal{Q}^g_+$, the condition $s_+ \cdot \psi = 0$ yields that $e_+ \cdot e_- \cdot w=0$ on $U$  which by multiplication with $e_-$ implies that $e_- \cdot w=0$. However, by (\ref{gdg}) and (\ref{h5}) it follows that on $U$ we have $D^g{\ph} = -n \cdot \widetilde{\Phi}^g ({proj}_-^g (\psi)) = 0$. Thus, $\ph$ is on $U$ both harmonic and a twistor spinor and therefore parallel wrt. $g$. Conversely, by the same argumentation every parallel spinor $\ph \in \Gamma(S^g)$ satisfies $s_+ \in \text{ker }{\psi}$. Moreover, one has in the language of Theorem \ref{gg} that $L=\text{pr}_{TM} (\text{ker }\psi \cap \mathcal{I}_-^{\bot})=\text{ker }\ph$ on $\widetilde{M}$ from which integrability of ker $\ph$ on $\widetilde{M}$ follows by Theorem \ref{gg}.
$\hfill \Box$ \\
\newline

\textbf{Proof of Theorem \ref{tss}: }
We first recall some algebraic facts for the relevant dimensions: In split signatures $(p,q) \in \{(m+1,m),(m,m) \}$ a real spinor $v \in \Delta^{\R}_{p,q}$ is called pure if dim ker $v = m$, i.e. ker $v=\{x \in \R^{p,q} \mid x \cdot v = 0 \}$ is maximal isotropic. Pure spinors do always exist, cf. \cite{har}, and they form a single orbit in $\Delta_{m+1,m}^{\R}$ or two orbits in $\Delta_{m,m}^{\R,+} \cup \Delta_{m,m}^{\R,-}$ under action of the spin group. It is known from \cite{br} that in signatures $(2,2)$, $(3,2)$ and $(3,3)$ every nonzero real (half)-spinor is pure, whereas in signatures $(4,4)$ and $(4,3)$ the orbit of pure spinors coincides with the set of nonzero (half-)spinors that have zero length wrt. $\langle \cdot, \cdot \rangle_{\Delta}$.\\
Let $\ph \in \text{ker }P^g$ be a nontrivial twistor half-spinor in signature $(2,2)$. The associated parallel spin tractor $\psi \in \mathcal{S}(M)$ is nowhere nonzero and therefore pointwise pure, i.e. dim ker $\psi$ = 3. Proposition \ref{tms} now yields the statement for signature $(2,2)$.\\
In signatures $(3,2)$ and $(3,3)$ the associated spin (half-)tractor $\psi$ to a nontrivial real twistor (half-)spinor $\ph$ is pure iff $\langle \psi, \psi \rangle_{\mathcal{S}} = 0$. However, \cite{hs} shows that this is equivalent to $\text{const.}=\langle \ph, D^g \rangle_{S^g} = 0$ for one (and hence for all) $g \in c$. That is, under the assumptions of Theorem \ref{tss} we have dim ker $\psi$=3 resp. 4 and again Proposition \ref{tms} applies.
$\hfill \Box$ \\
\newline

The parallel spinor fields arising via Theorem \ref{tss} are pointwise pure. \cite{kath,br} give a local normal form of the metric in this situation: For $(M,h)$ a pseudo-Riemannian spin manifold of split signature $(m+1,m)$ admitting a real pure parallel spinor field in $\Gamma(M,S^h)$, one can find for every point in $M$ local coordinates  $(x,y,z)$ , $x=(x_1,...,x_m)$, $y=(y^1,...,y^m)$ around this point such that 
\begin{align} h = -dz^2 - 4 \sum_{i=1}^m dx_i dy^i - 4 \sum_{i,j=1}^m g_{ij} dy^i dy^j, \label{purre1} \end{align}
where $g_{ij}$ are functions depending on $x,y$ and $z$ and satisfying
\begin{align}g_{ij} = g_{ji} \text{ for }i,j=1,...,m \text{,  }\sum_{i=1}^m \frac{\partial g_{ik}}{\partial x_i} = 0 \text{ for } k=1,...,m. \label{pure2} \end{align}
Conversely, if one uses (\ref{purre1}) and (\ref{pure2}) to define a metric $h$ on a connected open set $U \subset \R^{2m+1}$, then $(U,h)$ is spin and admits a real pure parallel spinor. $Hol(U,h)$ is contained in the image under the double covering $\lambda$ of the identity component of the stabilizer of a real pure spinor. $h$ is not necessarily Ricci-flat. Similar statements hold in case $(p,q)=(m,m)$, where one has to omit the last coordinate etc.\\
\newline
Consequently, real twistor spinors in signature $(3,2)$ fall in two disjoint classes distinguished by the constant $d:= \langle \ph, D^g \rangle_{S^g}$. For $d=0$, the local model is given by (\ref{purre1}) and the distribution ker $\ph$ is integrable, whereas for $d \neq 0$, the distribution ker $\ph$ is generic and the conformal structure can be recovered from it via a Fefferman construction, see \cite{hs1}.

\small
\bibliographystyle{amsalpha}
\bibliography{redhol}

\medskip 

\medskip

\textsc{Andree Lischewski\\
Humboldt-Universit\"at zu Berlin, Institut f\"ur Mathematik\\
Rudower Chaussee 25, Room 1.310, D12489 Berlin.\\
E-Mail: }\texttt{lischews@mathematik.hu-berlin.de}
\end{document}